\documentclass{amsart}
\usepackage{amsfonts}

\usepackage{amsmath}
\usepackage{amssymb}
\usepackage{amsthm}
\usepackage{verbatim}

\numberwithin{equation}{section}

\DeclareMathOperator{\re}{Re}

\begin{document}
\title[Coefficient invariances of Convex functions]{Coefficient invariances
of Convex functions}

\begin{abstract}
We summarise known sharp bounds for coefficient invariances for convex
functions, and suggest some significant open problems
\end{abstract}

\begin{abstract}
For convex univalent functions we give instances where the sharp bound for
various coefficient functionals are identical to those for the corresponding
bound for the inverse function . We give instances where the sharp bounds
differ and also suggest some significant open problems.
\end{abstract}

\title[coefficient invariances for convex functions]{coefficient invariances
for convex functions}
\author[D. K. Thomas]{Derek K. Thomas}
\address{Derek K. Thomas, Department of Mathematics, Swansea University Bay
Campus, Swansea, SA1 8EN, United Kingdom.}
\email{d.k.thomas@swansea.ac.uk}
\keywords{Univalent, Convex functions, Carath\'eodory functions, Inverse,
coefficients functions }
\subjclass{30C45}
\maketitle

%\maketitle

%\date{October 17, 2021}

\section{Introduction}

%\bigskip
Let ${\mathcal{A}}$ denote the class of analytic functions $f$ in the unit
disk $\mathbb{D}=\{ z\in\mathbb{C}: |z|<1 \}$ normalized by $%
f(0)=0=f^{\prime }(0)-1$. Then for $z\in\mathbb{D}$, $f\in {\mathcal{A}}$
has the following representation
\begin{equation*}  \label{A;01}
f(z) = z+ \sum_{n=2}^{\infty}a_n z^n.
\end{equation*}
\smallskip

Let ${\mathcal{S}}$ denote the subclass of all univalent (i.e., one-to-one)
functions in ${\mathcal{A}}$. \smallskip

The most significant subclasses of $\mathcal{S}$ are the classes $\mathcal{S}%
^*$ of Starlike and $\mathcal{C}$ of Convex functions, with $\mathcal{C}$
defined as follows. \bigskip

$f\in \mathcal{C}$ if, and only if, $f\in \mathcal{A}$ and

\bigskip

\begin{equation*}
\Re \left(1+\dfrac{zf^{\prime \prime }(z)}{f^{\prime }(z)}\right)>0,\quad
z\in \mathbb{D,}
\end{equation*}
\smallskip

A classical result using the Carath\'eodory functions shows that if $f\in
\mathcal{C}$, then $|a_n|\le 1$ for $n\ge2$, see e.g \cite{Pom}.

\section{Inverse coefficients}

Each $f\in\mathcal{S}$ possess an inverse function $f^{-1}$ given by

\begin{equation*}
f^{-1}(w)=w+\sum_{n=2}^{\infty}A_{n}w^n,
\end{equation*}
valid in some set $|w|<r_0(f)$ \bigskip

Although it took from 1916 to 1985 to solve the Bieberbach conjecture, in
1923 Lowner found the sharp upper bound for $|A_n|$ for all $n\ge2, $ by
showing that

\begin{equation*}
|A_n|\le \dfrac{\Gamma[2n+1]}{\Gamma[n+1]\Gamma[n+1]}.
\end{equation*}

Since the Koebe function $k(z)$ is also starlike, Lowner's bound is also
sharp for the class $\mathcal{S}^*$. \bigskip

However since the Koebe function does not belong to $\mathcal{C}$, other
sharp bounds must hold for $|A_n|$ when $f\in\mathcal{C}$, and a complete
solution to this problem appears difficult to find. \bigskip

The first indication that this was not a straightforward, but interesting
problem occurs in a paper in 1979 by Kirwan and Schober \cite{KS}, who
showed that there exists a function in $\mathcal{C}$ such that $|A_n|>1$ for
$n\ge10$. \bigskip

Then in 1982, Libera and Zlotkiewicz \cite{LZ} proved the following.

\bigskip %\begin{theorem}\label{Ye}
Let $f\in\mathcal{C}$ and have inverse function given by

\begin{equation*}  \label{Inverse}
f^{-1}(w)=w+\sum_{n=2}^{\infty}A_{n}w^n,
\end{equation*}

then for $2\le n\le7$, the following sharp inequalities hold.
\begin{equation*}
|A_n|\le1.
\end{equation*}

\bigskip

Subsequently in 1984, Campschroer \cite{Camp} showed that if
%\begin{theorem}\label{Ye}
$f\in\mathcal{C}$ and has inverse function given by (\ref{Inverse}), then $%
|A_8|\le1,$ and the inequality is sharp.

\bigskip

Thus the first questions that arise are, what is the sharp bound for $|A_9|$%
, and what are the sharp bounds for $|A_n|$ when $n\ge10$?

\bigskip

The next obvious question is to ask if there are any other invariance
properties amongst coefficient functionals in $a_n$ and $A_n$?

\smallskip

If it turns out that there are other invariant properties, then \textbf{why
is this so?} \bigskip

We will see that there are instances where sharp bounds can be found for
functionals concerning $A_n$ ($|A_n|$ for instance) which are different for
the sharp bound for the corresponding functional concerning $a_n$ ($|a_n|$
for instance)?

%{\center{Convex functions of order $\alpha$}}
\bigskip

Perhaps the most natural generalisation to the class $\mathcal{C}$ of convex
functions is the class of convex functions of order $\alpha$ defined as
follows. \bigskip

For $0\le\alpha<1$, denote by ${\mathcal{C}}(\alpha)$ the subclass of ${%
\mathcal{C}}$ consisting of convex functions of order $\alpha$ i.e., $f \in {%
\mathcal{C}}(\alpha)$ if, and only if, for $z\in\mathbb{D}$
\begin{equation*}
\re \left\{1+ \frac{zf^{\prime \prime }(z)}{f^{\prime }(z)} \right\} >\alpha.
\end{equation*}
\bigskip

\bigskip

The first obvious question to consider is to look for invariance amongst the
initial coefficients, where we at once encounter non-invariance between $%
|a_3|$ and $|A_3|$.

\medskip

Using elementary techniques, it is a simple exercise using well-known tools
to prove the follow inequalities, all of which are sharp

\medskip

If $f\in \mathcal{C}(\alpha),$ then\newline

\begin{equation*}
|a_{2}|, |A_{2}| \leq 1-\alpha, \quad \text{and} \quad |a_3|\le \dfrac{%
(3-2\alpha)(1-\alpha)}{3}.
\end{equation*}

\begin{equation*}
|A_{3}| \leq \left\{
\begin{array}{ll}
\dfrac{(3-4\alpha)(1-\alpha)}{3}, & \hbox{$0\le\alpha\leq\dfrac{1}{2}$,} \\
&  \\
\dfrac{1-\alpha}{3}, & \hbox{$\dfrac{1}{2}\leq\alpha<1$.}%
\end{array}
\right.
\end{equation*}

%{\center{Strongly convex functions}}
\bigskip

In 2016, Thomas and Verma \cite{TV} demonstrated some invariance properties
amongst the class of strongly convex functions defined as follows. \bigskip

%\center{Definition 1}
\bigskip

For $0<\beta\le1$, denote by ${\mathcal{C}}^{\beta}$ the subclass of ${%
\mathcal{C}}$ consisting of strongly convex functions i.e., $f \in {\mathcal{%
C}}^{\beta}$ if, and only if, for $z\in\mathbb{D}$
\begin{equation*}  \label{def1}
\Big|\arg \Big(1+ \frac{zf^{\prime \prime }(z)}{f^{\prime }(z)} \Big )|\le%
\dfrac{\beta \pi}{2}.
\end{equation*}
\bigskip

%{\center {Theorem 2}}
%\bigskip

We give first some invariant properties amongst the initial coefficients
proved for $f\in \mathcal{C}^{\beta}$ in \cite{TV}

\medskip

Let $f\in \mathcal{C}^{\beta}$ , then\newline
\begin{equation*}
|a_{2}|, |A_{2}| \leq \beta, \quad \quad \quad |a_{3}|,|A_{3}| \leq \left\{
\begin{array}{ll}
\dfrac{\beta}{3}, & \hbox{$0<\beta\leq\dfrac{1}{3}$,} \\
&  \\
{\beta}^2, & \hbox{$\dfrac{1}{3}\leq\beta\leq1$.}%
\end{array}
\right.
\end{equation*}%
\newline

\begin{equation*}
|a_{4}|,|A_{4}| \leq \left\{
\begin{array}{ll}
\dfrac{\beta}{6}, & \hbox{$0<\beta\leq\sqrt{\dfrac{2}{17}}$,} \\
&  \\
\dfrac{\beta}{18}(1+17{\beta}^2), & \hbox{$\sqrt{\dfrac{2}{17}}\leq\beta%
\leq1$.}%
\end{array}
\right.
\end{equation*}

All the inequalities are sharp.

%{\center{A Fekete-Szeg\"o Invariance }}
\bigskip

If $f\in \mathcal{C}^{\beta}$, then for any complex number $\nu$,\newline
\begin{equation*}
\left|a_{3}-\nu{a_{2}^{2}}\right|, \left|A_{3}-\nu{A_{2}^{2}}%
\right|\leq\max\left\{\frac{\beta}{3}, \ {\beta}^2\left|1-\nu\right|\right\}.
\end{equation*}

Both inequalities are sharp. \bigskip

%{\center {Logarithmic Coefficients}}
%\bigskip

For $f\in\mathcal{A}$, the logarithmic coefficients $\gamma_n$ of $f(z)$ are
defined by

\begin{equation*}
\log \dfrac{f(z)}{z}=2\sum_{n=1}^{\infty}\gamma_nz^n.
\end{equation*}

They play a central role in the theory of univalent functions, and formed
the basis of de Brange's proof of the Bieberbach conjecture. \medskip

We make the following definition. \medskip

Let $f\in \mathcal{C}^{\beta}$, and $\log \dfrac{f^{-1}(\omega)}{\omega}$ be
given by
\begin{equation*}
\log \dfrac{f^{-1}(\omega)}{\omega}=2\sum_{n=1}^{\infty} c_{n}{\omega}^{n}.
\end{equation*}

%{\center{Theorem 3}}
%\bigskip

It was also shown in \cite{TV} that if $f\in \mathcal{C}^{\beta}$, then
%\bigskip

\begin{equation*}
|\gamma_1|, |{c}_{1}| \leq \dfrac{\beta}{2}, \quad \quad \quad |\gamma_2|, |{%
c}_{2}| \leq \left\{
\begin{array}{ll}
\dfrac{\beta}{6}, & \hbox{$0<\beta\leq\dfrac{2}{3}$,} \\
&  \\
\dfrac{{\beta}^2}{4}, & \hbox{$\dfrac{2}{3}\leq\beta\leq1$,}%
\end{array}
\right.
\end{equation*}%
\newline

\begin{equation*}
|\gamma_3|, |{c}_{3}| \leq \left\{
\begin{array}{ll}
\dfrac{\beta}{12}, & \hbox{$0<\beta\leq\sqrt{\dfrac{2}{5}}$,} \\
&  \\
\dfrac{\beta}{36}(1+5{\beta}^2), & \hbox{$\sqrt{\dfrac{2}{5}}\le
\beta\leq1$.}%
\end{array}
\right.
\end{equation*}
All the inequalities are sharp.

%{\center {Hankel Determinants}}
\bigskip

Clearly the more complicated the coefficient functional, the more difficult
will be the analysis, and finding invariance. \bigskip

We consider next the second Hankel determinants $H(2,2)(f)$ and $%
H(2,2)(f^{-1})$, defined by

\begin{equation*}
H(2,2)(f)=a_2 a_4-a_3^2,
\end{equation*}
{and}
\begin{equation*}
H(2,2)(f^{-1})=A_2 A_4-A_3^2
\end{equation*}
\bigskip

%{\center {Hankel Determinant Invariance}}
%\medskip

It was further shown by Thomas and Verma \cite{TV} that if $f\in\mathcal{C}%
^{\beta}$, then

\begin{equation*}
|H(2,2)(f)|,|H(2,2)(f^{-1})| \le
\begin{cases}
\dfrac{\beta^2}{9}, & \mbox \ 0<\beta\leq \dfrac{1}{3}, \\
&  \\
\dfrac{\beta(1+\beta)(1+17\beta)}{72(3+\beta)}, & \mbox\ \dfrac{1}{3}\le
\beta\le 1 \\
&
\end{cases}%
\end{equation*}

\noindent who claimed that all the inequalities are sharp. \medskip

Although the proofs of the positive results are correct, the claim that the
second inequality is sharp is false.

Note however that when $\beta=1,$ i.e., $f\in\mathcal{C}$, $%
|H(2,2)(f)|,|H(2,2)(f^{-1})|\le\dfrac{1}{8},$ and these inequalities are
sharp.

%{\center {Hankel Determinant Invariance (Correction)}}
\bigskip

Although the methods used in the proofs of the above invariance are correct,
they are not strong enough to give sharp bounds for the second inequality.
The following correction was subsequently given by Lecko, Sim and Thomas
\cite{LST}. \medskip

If $f\in\mathcal{C}^{\beta}$, then the following sharp inequalities hold.

\begin{equation*}
|H(2,2)(f)|,|H(2,2)(f^{-1})| \le
\begin{cases}
\dfrac{\beta^2}{9}, & \mbox \ 0<\beta\leq \dfrac{1}{3}, \\
&  \\
\dfrac{\beta^2(1+\beta)(17+\beta)}{72(2+3\beta-\beta^2)}, & \mbox\ \dfrac{1}{%
3}\le \beta\le 1. \\
&
\end{cases}%
\end{equation*}

\medskip

Note again that when $\beta=1,$ i.e., $f\in\mathcal{C}$, $%
|H(2,2)(f)|,|H(2,2)(f^{-1})|\le\dfrac{1}{8.}$

%{\center { Successive Coefficient Differences}}
\medskip

The problem of finding sharp bounds for the difference of successive
coefficients $|a_{n+1}|-|a_n|$ for functions in $\mathcal{S}$ represents one
of the most difficult areas of study in univalent function theory, with the
only sharp bound known so far is when $n=2$, where Duren prove the rather
curious sharp inequality

\begin{equation*}
-1 \leq |a_3| - |a_2| \leq \frac{3}{4} + e^{-\lambda_0}(2e^{-\lambda_0}-1) =
1.029\cdots,
\end{equation*}
where $\lambda_0$ is the unique value of $\lambda$ in $0 < \lambda <1$,
satisfying the equation $4\lambda = e^{\lambda}$.

\bigskip

Although Leung found the complete solution when $f$ is starlike by showing
that $||a_{n+1}|-|a_n||\le1$, the problem for convex functions remains
mostly open.

%{\center { Successive Coefficient Differences for $f\in\mathcal{C}$}}
\medskip

In 2016, Ming and Sugawa \cite{Li} made the first advance in this problem by
finding the following sharp bounds when $n\ge2$
\begin{equation*}
|a_{n+1}|-|a_n|\le \dfrac{1}{n+1},
\end{equation*}
and also proved the sharp lower bounds

\begin{equation*}
|a_{3}|-|a_2|\ge -\dfrac{1}{2},\quad \text{and}\quad |a_{4}|-|a_3|\ge -%
\dfrac{1}{3}
\end{equation*}

Thus in particular we have the sharp bounds

\begin{equation*}
-\dfrac{1}{2}\le|a_{3}|-|a_2|\le \dfrac{1}{3} \quad \text{and}\quad -\dfrac{1%
}{3}\le|a_{4}|-|a_3|\le \dfrac{1}{4}.
\end{equation*}

%{\center { Inverse  Coefficient Differences for $f\in\mathcal{C}$}}
\medskip

Sim and Thomas in 2020 \cite{ST} proved the following inequalities hold for
the inverse coefficients, thus demonstrating another example of invariance.
\begin{equation*}
-\dfrac{1}{2}\le|A_{3}|-|A_2|\le \dfrac{1}{3},
\end{equation*}

\bigskip

\noindent and have recently shown, (the proof of which requires much more
complicated analysis), that the following sharp inequalities hold

\begin{equation*}
-\dfrac{1}{3}\le|A_{4}|-|A_3|\le \dfrac{1}{4}.
\end{equation*}
\medskip

Any advances when $n\ge 4$ would require deeper methods of proof.

\medskip

It is clear from the above, that the more complicated the functional and
class of convex functions considered, the less likely is it that there will
be invariance. Also functionals containing coefficients $a_n$ and $A_n$ for $%
n\ge3$ will similarly be difficult to deal with. \medskip

Sim and Thomas (\cite{ST2}, Proposition 1) have recently given a general
lemma concerning functions of positive real part, which provides a tool
enabling coefficient differences to be found when $n=2$, which can be
applied to many subclasses of univalent functions . This can also be used to
consider coefficient differences of the inverse function. However the lemma
only applies when $n=2$.

\medskip

%{\center { Other Invariances }}
\medskip

In recent years it has become fashionable to consider subclasses of convex
(and starlike) functions where the function $p(z)$ is specified, often
having a range with some interesting geometrical property. \medskip

Examples of recently discussed subclasses of convex functions, where some
initial invariance properties have been found are as follows.

%{\center { Convex functions related to the Exponential function (1) }}
\medskip

A most natural class of convex functions related to the exponential function
is the class $\mathcal{C}_{E}$ defined as follows

\begin{equation*}
\mathcal{C}_{E}=\left\{ f\in \mathcal{A}:1+\frac{zf^{\prime \prime }(z)}{%
f^{\prime }(z)}\prec e^{z}\right\}.
\end{equation*}
\medskip

Some initial coefficients results for the class $\mathcal{C}_{E}$ were given
by Zaprawa \cite{Zap}, and similar analysis shows that the following
invariance properties hold.

\begin{equation*}
|a_{2}|, |A_2|\le \dfrac{1}{2} \quad \text{and}\quad |a_{3}|, |A_3|\le
\dfrac{1}{4}\quad \text{and}\quad |a_{4}|, |A_4|\le \dfrac{17}{144}.
\end{equation*}

All the inequalities are sharp.

%{\center { More invariances for  $f\in\mathcal{C}_{E}$}}
\medskip

Using well-known methods it is also possible to prove the following
inequalities, both of which are sharp.

\begin{equation*}
|a_2 a_3-a_4|, |A_{2} A_{3}-A_4|\le \dfrac{1}{12}.
\end{equation*}

\medskip

Similarly the following sharp invariances hold for the second Hankel
determinants for $f\in\mathcal{C}_{E}$ \cite{Zap},

\begin{equation*}
|H(2,2)(f)|, |H(2,2)(f^{-1})|\le \dfrac{73}{2592}.
\end{equation*}
\medskip

%{\center { Convex functions related to the Exponential function (2) }}
\medskip

A class $\mathcal{C}_{SG}$ of convex functions which exhibits some
interesting invariance properties is related to a modified sigmoid function
and is defined by
\begin{equation*}
\mathcal{C}_{SG}=\left\{ f\in \mathcal{A}:1+\frac{zf^{\prime \prime }(z)}{%
f^{\prime }(z)}\prec \frac{2}{1+e^{-z}}\right\} .
\end{equation*}

Here the function $2/(1+e^{-z})$ is a modified sigmoid function which maps $%
\mathbb{D}$ onto the domain $\Delta _{SG}=\left\{ w\in \mathbb{C}%
\right\}:\left\vert \log \left( w/\left( 2-w\right) \right) \right\vert <1\}
.$

\medskip

The class $\mathcal{C}_{SG}$ was first discussed in \cite{DKT} where some
invariance properties were found. In particular the following invariance
properties hold.

\begin{equation*}
|a_{2}|, |A_2|\le \dfrac{1}{4} \quad \text{and}\quad |a_{3}|, |A_3|\le
\dfrac{1}{12}\quad \text{and}\quad |a_{4}|, |A_4|\le \dfrac{1}{24}.
\end{equation*}

All the inequalities are sharp.

%{\center { Further invariance properties for functions in $\mathcal{C}_{SG}$}}
\medskip

Further invariance properties for functions in $\mathcal{C}_{SG}$ also hold
(see \cite{DKT1}), where the following inequality for $|a_2 a_3-a_4|$ was
proved, and the inequality for $|A_2 A_3-A_4|$ follows using similar
methods, \medskip

\begin{equation*}
|a_{2}a_3-a_4|, |A_2 A_3-A_4|\le \dfrac{1}{24}.
\end{equation*}

\noindent Both inequalities are sharp. \medskip

As already mentioned, finding sharp bounds for the differences of
coefficients can present difficulties and next we give an example of proved
non-invariance for $\mathcal{C}_{SG}$, noting that \textbf{all the
inequalities are sharp}.

\begin{equation*}
-\dfrac{5}{24}\le |a_{3}|-|a_2|\le \dfrac{1}{12},\quad \text{and} \quad -%
\dfrac{1}{4}\le |A_{3}|-|A_2|\le \dfrac{1}{12}.
\end{equation*}
\medskip

\medskip

As mentioned above other choices of $p(z)$ have been made primarily to
define some kind of interesting geometric property of the range of $p(z)$,
for example

\smallskip

(i) \ $p(z)=1+\dfrac{4}{3}z+\dfrac{2}{3}z^2$ gives a cardioid domain.
\smallskip

(ii) $p(z)=1+\dfrac{4}{5}z+\dfrac{1}{5}z^4$ gives a 3-leaf petal shaped
domain.

\medskip

It is very likely that similar invariances hold for some functionals in
these classes.

\medskip

Finally we note a recent interesting example of non-invariance. \bigskip

It was shown by Sim, Zaprawa and Thomas in 2021 \cite{SZT} that for $f\in
\mathcal{C}(\alpha),$

\begin{equation*}
|H(2,2)(f)|\le \dfrac{(1-\alpha)^2(6+5\alpha)}{48(1+\alpha)},
\end{equation*}

\noindent and that this inequality is sharp, thus solving a long-standing
problem. \smallskip

\noindent In the same paper it was shown that for $f\in \mathcal{C}(\alpha),$

\begin{equation*}
|H(2,2)(f^{-1})| \le
\begin{cases}
\dfrac{1}{96}(12-28\alpha+19\alpha^2), & \mbox \ \alpha\in [0,\dfrac{2}{5}],
\\
&  \\
\dfrac{1}{9}{(1-\alpha)^2}, & \mbox\ \alpha\in [\dfrac{2}{5},\dfrac{4}{5}],
\\
&  \\
\dfrac{\alpha(1-\alpha)^2(19\alpha-8)}{48(1+\alpha)(2\alpha-1)}, & \mbox\ %
\alpha\in [\dfrac{4}{5},1].%
\end{cases}%
\end{equation*}
\smallskip

\noindent and that all the inequalities are sharp.

\bigskip

\noindent Thus unless $\alpha=0$, we have non-invariance.

\bigskip

We note here that when $f\in \mathcal{C}^(\beta),$ there is invariance
between $|H(2,2)(f)|$ and $|H(2,2)(f^{-1})|$ for all $\beta\in (0,1]$, which
is curious since the definition of the class $\mathcal{C}^(\beta)$ involves
the power $p(z)^{\beta}.$

\bigskip

Checking for invariances is a matter of applying available well-known tools
to functionals.

\smallskip

But the \textbf{real problem} is to discover \textbf{WHY} invariances occur
in the classes and functionals considered.

\smallskip

The answer is probably that invariances are both class and functional
dependant, and that there is no simple rule.

\bigskip

Recall that in 2016, Ming and Sugawa \cite{Li} proved that if $f\in\mathcal{C%
}$, then when $n\ge2$
\begin{equation*}
|a_{n+1}|-|a_n|\le \dfrac{1}{n+1},
\end{equation*}
and that the bounds are sharp. \bigskip

In the same paper Ming and Sugawa \cite{Li} further proved that when $n\ge4$

\begin{equation*}
-\dfrac{1}{n}<|a_{n+1}|-|a_n|\le \dfrac{1}{n+1},
\end{equation*}

\smallskip

\noindent thus

\begin{equation*}
||a_{n+1}|-|a_n||=\mathcal{O} \big(\dfrac{1}{n}\big), \ \ \text{as}\ n\to
\infty.
\end{equation*}

\smallskip

Is it true that if $f\in \mathcal{C}$, then

\begin{equation*}
||A_{n+1}|-|A_n||=\mathcal{O} \big(\dfrac{1}{n}\big), \ \ \text{as}\ n\to
\infty?
\end{equation*}
\bigskip

\smallskip

Next note that Kowalczyk, Lecko and Sim \cite{KLS} have shown that if $f\in
\mathcal{C}$, then the third Hankel determinant
\begin{equation*}
|H_{3,1}(f)| =|a_3(a_2a_4 - a_3^2) - a_4(a_4 - a_2a_3) + a_5(a_3 -
a_2^2)|\le \dfrac{4}{135},
\end{equation*}

\noindent and that this inequality is sharp. \smallskip

It is true that the following sharp inequality holds?

\begin{equation*}
|H_{3,1}(f^{-1})| =|A_3(A_2A_4 - A_3^2) - A_4(A_4 - A_2A_3) + A_5(A_3 -
A_2^2)|\le \dfrac{4}{135}.
\end{equation*}

\end{document}